\documentclass[12pt,leqeq]{article}

\usepackage{amssymb,amsthm}
\usepackage[centertags]{amsmath}
\usepackage{graphics}

\usepackage{amsmath}

\usepackage{amscd}

\usepackage{epsf,graphicx}

\usepackage[all]{xy}

\usepackage{amsfonts}

\usepackage{newlfont}

\usepackage{feynmp}

\newtheorem{thm}{Theorem}

\newtheorem{lem}[thm]{Lemma}

\newtheorem{cor}[thm]{Corollary}

\newtheorem{prop}[thm]{Proposition} 

\newtheorem{defn}[thm]{Definition}

\newtheorem{exmp}[thm]{Example}

 \newcommand{\Z}{\mathbb{Z}}

 \newcommand{\N}{\mathbb{N}}

\begin{document}

\setlength{\baselineskip}{16pt}

\title{$A^{N}_{\infty}$-algebras}

\author{Mauricio Angel and Rafael D\'\i az}

\maketitle


\begin{abstract}
We study higher depth algebras. We introduce several examples of such structures
starting from the notion of $N$-differential graded algebras and build up to
the concept of $A_{\infty}^N$-algebras.
\end{abstract}

\section{Introduction}

In this paper we continue our program aim to study
algebraic structures of higher depth, i.e., algebraic structures such that
the axioms satisfied by the structural operations on it involve triple
or higher order compositions. The simplest and most prominent algebraic
structure of such type are $N$-complexes. According to Mayer \cite{Ma}  and Kapranov
\cite{Kap} a $N$-complex is a graded vector space $V$ together with a
degree one map $d:V \rightarrow V$ such that $d^N=0.$ This condition
involves $N$ compositions of the structural map $d$, this makes
the definition of $N$-complexes quite different from that
of a complex which involves only two compositions of
the operator $d$, and so it is quadratic, as most algebraic structures we are familiar with.
In this paper we introduce other examples of algebraic structures
of higher depth. We put special interest in the operadic description of such structures and
in the corresponding deformation theory.\\

Algebraic structures of  higher depth arise when one
defines an associative product on a $N$-complex. In this situation there are two consistent
choices. One choice begins fixing $q$ a primitive $N$-root of unity. Then one defines the
notion of $q$-differential graded algebras, which are associative graded algebras provided
with a $N$-differential satisfying a twisted graded Leibnitz rule. The condition
of associativity and the Leibnitz rule are quadratic, but the condition for the $N$-differential is of
higher depth. This sort of algebraic structure was first discussed by Kapranov and Dubois-Violette
in \cite{DV1}, \cite{DV2}, \cite{Kap}, and has been furhter studied in several works, among them
\cite{AK}, \cite{DV3}, \cite{DVK}, \cite{KN}, \cite{KaWa}.
Deformations of $q$-differential graded algebras are controlled
by the $q$-analogue of the Maurer-Cartan equation, which was defined for general $N$ and
explicitly computed for small values of $N$ in \cite{AnDi3}. \\

The other choice introduced in \cite{AnDi1} is that of $N$-differential graded algebras, which are $N$-complexes
provided with an associative product such that the $N$-differential satisfies the graded Leibnitz
rule. Again this is a sort of hybrid algebraic structure with the higher depth condition coming
from the $N$-differential. We remark that by now there are quite a few known
examples of algebras of this sort, we summarize these examples in Section \ref{ndgla}.
In this work we introduce the graded operad $N$-$dga$ whose
algebras are $N$-differential graded algebras and compute its generating series. We study a couple
of related algebraic structures, namely, $N$-differential graded Lie algebras and $N$-codifferential graded
coalgebras.  We introduce the graded operad $N$-$dgla$ whose algebras are $N$-differential graded Lie algebras
and compute its generating function. In Section \ref{co} we study deformations of $N$-codifferential graded coalgebras
into $M$-codifferential graded coalgebras
and show that they are controlled by the $(N,M)$ Maurer-Cartan equation.\\

All algebraic structures mention so far
are of higher depth only because of the condition imposed on a $N$-differential. Our
next goal in this paper is to introduce new families of algebras of higher depth, which are homogeneous
in the sense that in each axiom present in their definitions,
exactly $N$ structural operations are involved. In order to motivate our definitions
it is convenient to switch to the world of graded coalgebras and  codifferentials. Recall
that a differential graded algebra structure on a graded vector space $A$
is the same as a coderivation $\delta$ on the graded coalgebra $T^{\leq 2}(A[1])$
satisfying $\delta^2=0$. If we relax the latter condition and demand instead that
$\delta^N=0$, then we arrive to the notion of differential graded algebra of depth $N$.
Any such coderivation is determined by a couple of maps $d:A[1] \rightarrow A[1]$ and
$m:A[1] \otimes A[1] \rightarrow A[1]$. If the product $m$ vanishes
one recovers the definition of a $N$-complex. If instead $d$ is zero then we arrived to the
notion of a $N$-associative algebra.\\

We construct the operad whose algebras are differential graded algebras of depth $N$,
and also the operad whose algebras are $N$-associative algebras. We show that infinitesimal
deformations of $N$-associative algebras into $M$-associative algebras, $M \geq N$, are controlled
by a well-defined cohomology group. We write explicitly the conditions determining
if an associative algebra admits a non-trivial deformation into a $3$-associative algebra.
In the final section we introduce the notion of $A_{\infty}^{N}$-algebras which generalizes both
$N$-associative algebras and $A_{\infty}$-algebras.
$A_{\infty}^{N}$-algebras are defined as nilpotent coderivations on the full cotensor coalgebra
$T(A[1])$. We introduce the operad for $A_{\infty}^{N}$-algebras, study infinitesimal
deformations of $A_{\infty}^{N}$-algebras, and show that the moduli space of deformations of $A_{\infty}^{N}$-algebras
into $A_{\infty}^{M}$-algebras is controlled by the $(N,M)$ Maurer-Cartan equation mentioned
above.\\

We close the introduction emphasizing that we are only beginning the study
higher depth algebras. Much more work, both theoretical and practical, is needed
in order to haver a better grasp of the meaning and applications of such structures.
We believe that new forms of infinitesimal symmetries will be uncover along this
line of thought through the notion of Lie $N$-algebroids to be developed in
\cite{ACD}.

\section{N-differential graded algebras}\label{ndgla}

Throughout this paper we  work with the category $gvect$ of
graded $k$-vector spaces $V$ over a field $k$ of characteristic
zero. The degree of an homogeneous element $v \in V$ is denote by
$\bar{v}\in\Z$. $V[k]$ denotes the $\Z$-graded vector space such
that $(V[k])^i=V^{i+k}$ for  $i\in\Z$. The superdimension of a graded vector space $V$ is given by
$$sdim(V)=\sum_{i \in \mathbb{Z}}dim(V_{2i}) - \sum_{i \in \mathbb{Z}}dim(V_{2i+1}).$$

\begin{defn}
A $N$-complex is a graded vector space $V$ together with a degree one map
$d:V \longrightarrow V$ such that $d^N=0.$ An $N$-complex is said to be a proper
$N$-complex if $d^{N-1}\neq 0$.
\end{defn}

\begin{defn}\label{Ndga}
 A   $N$-differential graded algebra (N-dga) is
a triple $(A,m,d)$, where $A$ is a graded vector space, $m:A \otimes A \to A$
and $d:A
\to A$  are maps of degree zero and one,
respectively, such that
\begin{enumerate}
\item $(A,m)$ is a graded associative algebra.

\item $d$ satisfies the graded Leibnitz rule
$d(ab)=d(a)b+(-1)^{\bar{a}}ad(b)$ for $a,b \in A.$

\item $d^{N}=0$.
\end{enumerate}
\end{defn}

We said that a $A$ is a proper $N$-dga if $d^{N-1}\neq0$.

\begin{defn}
\begin{enumerate}
\item Let $N^{il}$-dgvect be the category whose objects are $N$-differential graded
vector spaces for some $N$. Morphisms in $N^{il}$-dgvect are maps $T:V \rightarrow W$ such
that $dT=Td$.

\item Let $N^{il}$-dga be the category whose objects are $N$-differential graded
algebras for some $N$. Morphisms in $N^{il}$-dga are maps $T:V \rightarrow W$ such
that $dT=Td$ and $m(T \otimes T)=Tm$.

\end{enumerate}
\end{defn}

\smallskip

The following result justifies from the point of view of category theory our definition
of $N$-differential graded algebras.

\begin{lem}
Tensor product gives $N^{il}$-dgvect the structure of a monoidal category.
$N^{il}$-dga is the category of monoids in $N^{il}$-dgvect.
\end{lem}

We remark that by now there are several known examples of
nil-differential graded algebras, see \cite{ACD}, \cite{AnDi1}, \cite{AnDi2} and \cite{AnDi3}, which may be classified as
follows:

\begin{enumerate}

\item Differential graded algebras are the same as $2$-differential graded
algebras.

\item If $C$ is a $N$-complex, then $End(C)$ is a $(2n+1)$-differential graded algebra. There are plenty of examples of
$N$-complexes \cite{AK}, \cite{DV3},\cite{DVK}, \cite{Kap}, \cite{KaWa}, \cite{KN}.

\item $N$-flat connections and $N$-flat Riemannian metrics gives
rise to nil-differential graded algebras. A connection $A$ on a vector
bundle is called $N$-flat if and only if its curvature $F_A$ satisfies $(F_A)^N=0.$ \
A Riemaniann metric $g$ on a manifold is $N$-flat if the Levi-Civita connection on the tangent bundle is $N$-flat.

\item Differential forms of depth $N$ on affine manifolds of
dimension $m$ are $m(N-1)+1$ differential graded algebras.

\item One can generalize Sullivan's notion of algebraic differential forms on a simplicial sets,
to the notion of differential forms of depth $N$.  For finitely generated
simplicial sets these algebras are examples of Nil-differential graded algebras.

\item Zeilberger in \cite{Z} defines the notion of difference forms on the integral lattice in Euclidean space.
His construction can be generalized to define the algebra of difference forms of depth $N$ on
simplicial sets. Difference forms of depth $N$ on finitely generated simplicial sets
are twisted Nil-differential graded algebras.

\item Lie $N$-algebroids are examples of $N$-differential graded algebras.
\end{enumerate}

We are going to use the language of operads, the reader may consult
\cite{L} for more on operads. The generating series of an operad $O$ in $vect$
is $\sum_{n=0}^{\infty}dim(O(n))\frac{N^n}{n!}$, the generating series of an operad $O$ in $gvect$
is $\sum_{n=0}^{\infty}sdim(O(n))\frac{N^n}{n!}.$

\begin{defn}
Let $N\mbox{-}dga$ be the free graded operad generated by a degree one element
$d \in N\mbox{-}dga(1)$, and a degree zero  element $m \in N\mbox{-}dga(2)$,
subject to the relations $m^2=0$ and $dm=dm$, $d^N=0.$
\end{defn}

Graphically $N\mbox{-}dga$ is the free operad generated by the tree

\begin{figure}[h!]
\begin{center}
\includegraphics[height=1.5cm]{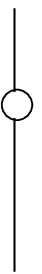}\label{a}
\end{center}
\end{figure}

\newpage
representing the $N$-differential and the tree

\begin{figure}[h!]
\begin{center}
\includegraphics[height=1.5cm]{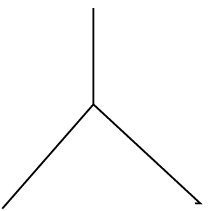}\label{b}
\end{center}
\end{figure}

representing the product, subject to the relations

\begin{figure}[h!]
\begin{center}
\includegraphics[height=1.5cm]{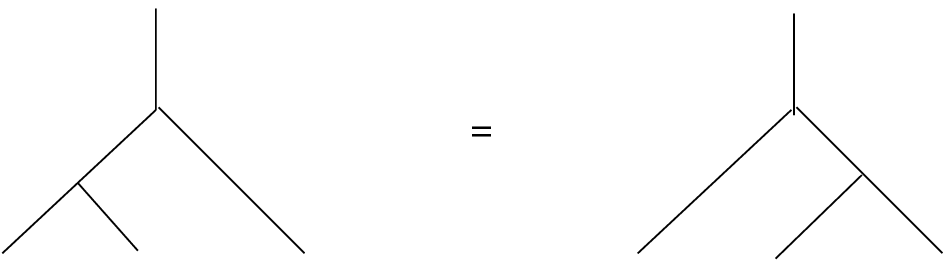}\label{dib4}
\end{center}
\end{figure}

\begin{figure}[h!]
\begin{center}
\includegraphics[height=1.5cm]{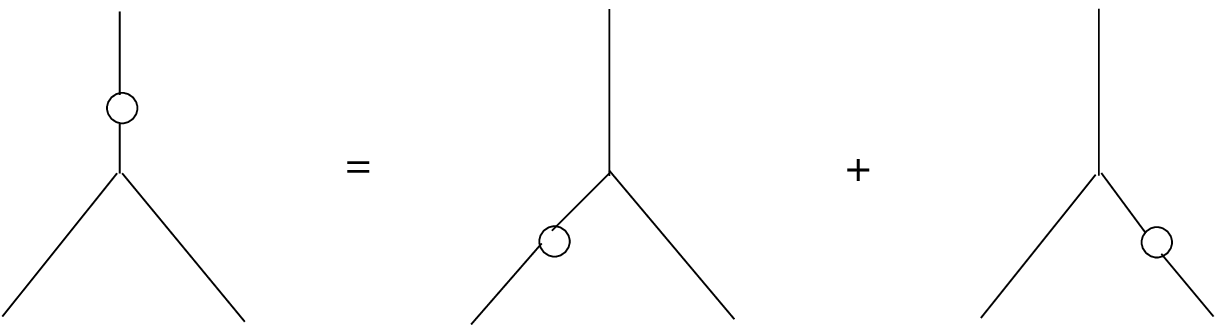}\label{dib3}
\end{center}
\end{figure}

\begin{figure}[h!]
\begin{center}
\includegraphics[height=2.5cm]{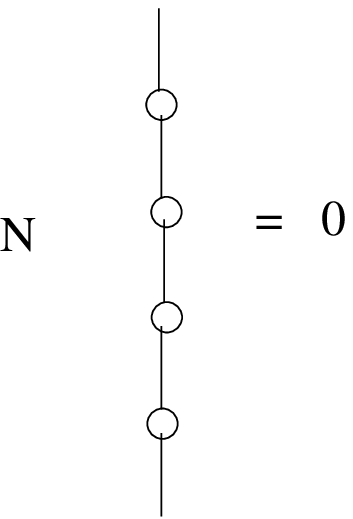}\label{dib5}
\end{center}
\end{figure}

The following result is proved using the graphical description above.

\begin{prop}
\begin{enumerate}
\item Algebras in $gvect$ over $N\mbox{-}dga$ are $N$-differential graded algebras.

\item $dim(N\mbox{-}dga(n))= n!N^{n}$. The generating series of $N\mbox{-}dga$ is $\frac{Nx}{1-Nx}$.

\item If $N$ is even then $sdim(N\mbox{-}dga(n))=0,$ for $n \geq 2$. The generating series of
$N\mbox{-}dga$ as a graded linear operad is $x$.

\item If $N$ is odd then $sdim(N\mbox{-}dga(n))=n!$, for $n \geq 1$. The generating series of
$N\mbox{-}dga$ as a graded operad is $\frac{x}{1-x}$.

\end{enumerate}
\end{prop}

Let us introduce the corresponding notion in the context of Lie algebras.

\begin{defn}
A ${N}$-differential graded Lie algebra ($N$-dgla) is a
$\Z$-graded vector space $L$ together with  a degree zero map $[\ ,\
]:L \otimes L \to L$ and a degree one map $d:L \to L$ such that
\begin{enumerate}
\item $(L,[\ ,\ ])$ is a graded Lie algebra.
\item $d[a,b]=[d(a),b]+(-1)^{\bar{a}}[a,d(b)]$.
\item $d^N=0$.
\end{enumerate}
\end{defn}

Let $\tau:\{1,2 \} \rightarrow \{1,2 \}$ be the non-trivial permutation on $\{1,2 \}.$

\smallskip

\begin{defn}
 Let $N\mbox{-}dgla$ be the free graded operad generated by
$d \in N\mbox{-}dgla(1)$ of degree one and $m \in N\mbox{-}dgla(2)$ of degree zero,
subject to the relations $m\tau=-m$, $m^2=0$, $dm=dm$ and $d^N=0$.
\end{defn}

\smallskip

\begin{prop}
\begin{enumerate}
\item Algebras over $N\mbox{-}dgla$ are $N$-differential graded Lie algebras.

\item $dim(N\mbox{-}dgla(n))=(n-1)!N^n$. The generating series of $N\mbox{-}dgla$ as a
linear operad is $ln(\frac{1}{1-Nx})$.

\item If $N$ is even the $superdim(N\mbox{-}dgla(n))=0,$ for $n \geq 2$. The generating series of
$N\mbox{-}dgla$ as a graded  operad is $x$.

\item If $N$ is odd the $superdim(N\mbox{-}dgla(n))=(n-1)!,$ for $n \geq 2$. The generating series of
$N\mbox{-}dgla$ as a graded  operad is $ln(\frac{1}{1-x})$.

\end{enumerate}
\end{prop}

\section{Deformation of coalgebras}\label{co}

In this section we introduce the notion of $N$-codifferential graded coalgebras,
and study the deformation theory of $N$-codifferentials.

\begin{defn}
A graded coalgebra is a $\Z$-graded vector space $C$
together with a degree zero map $\Delta:C\to C\otimes C$
such that the diagram
\[\xymatrix{
C  \ar[r]^{\Delta}\ar[d]_{\Delta} &C\otimes C\ar[d]^{\Delta\otimes 1} \\
C\otimes C  \ar[r]_{1\otimes\Delta} & C\otimes C\otimes C } \]
commutes. If in addition we have a degree zero map $\epsilon:C \to k$
such that the diagram
\[ \xymatrix{
 & C\ar[d]^{\Delta} & \\
C  \ar[ur]^{1}& C \otimes C\ar[l]^{\epsilon\otimes 1}
\ar[r]_{1\otimes\epsilon} & C \ar[ul]_{1} } \] commutes, then we
say that $C$ is a graded coalgebra with counit.
\end{defn}
\smallskip

\begin{defn}
A coderivation $\delta:C \to C$ on a graded coalgebra
$(C,\Delta)$ is a linear map such that the  diagram
\[\xymatrix{ C  \ar[r]^{\delta} \ar[d]_{\Delta} & C \ar[d]^{\Delta} \\
C \otimes C \ar[r]_{1\otimes\delta+\delta\otimes 1} & C \otimes C
}
\] commutes.  We denote by $Coder(C)$ the space of coderivations
on $C$.
\end{defn}

\begin{defn}
A $N$-codifferential graded coalgebra ($N$-cgc) is a pair
$(C,\delta)$ where $C$ is a $\Z$-graded coalgebra and $\delta:C
\to C$ is a degree one coderivation  such that $\delta^N=0$.
\end{defn}

A $1$-cgc is a graded coassociative coalgebra. A $2$-cgc is a
codifferential graded coalgebra. Below we use the fact that
a $N$-differential graded algebra $A$ may be regarded as a $N$
differential graded Lie algebra with bracket
$[a,b]=ab-(-1)^{\overline{a}\overline{b}}ba,$ for $a,b
\in A$.

\begin{prop}
Let $(C,\delta)$  be a $N$-cgc.
\begin{enumerate}
\item $End(C)$ is a $(2N-1)$-dga.
\item $Coder(C)$ is a $(2N-1)$-differential
graded Lie algebra.
\end{enumerate}
\end{prop}
\begin{proof} 1. Composition gives  $End(C)$  the structure of an
associative algebra. Differential $d$ on $End(C)$ is defined
by $d(f)=\delta\circ f-(-1)^{\bar{f}}f\circ\delta$, for $f \in End(C)$. $d$ and $m$ satisfy the graded Leibniz rule.
The identity
\[d^n(f)=\sum_{k=0}^n (-1)^{n-k}\binom{n}{k} \delta^k\circ f\circ \delta^{n-k}\]
implies that $(End(C),m,d)$ is a $(2N-1)$-dga.

\smallskip

2. By the previous remark $End(C)$ is a $(2N-1)$-dgla. It is easy
to show that $Coder(C)\subseteq End(C)$ is closed under Lie
brackets, thus it is a $(2N-1)$-dgla.
\end{proof}

\bigskip

We consider deformations of $N$-codifferentials and show
that such deformations are controlled by the $(N,M)$
Maurer-Cartan equation. Let  Artin be the category whose objects are
local $k$-algebras.
\begin{defn}\label{def}
Let $C$ be a $N$-cgc and $a$ and object in {\em $\mbox{Artin}$}
with maximal ideal $a_+$. A {\em $M$-deformation} of $C$ over $a$
is a $M$-cgc $C_a$ over $a$, with $M \geq N$, such that
$C_a/a_+C_a$ is isomorphic to $C$ as $N$-cgc.
\end{defn}

As usual in deformation theory the core of Definition \ref{def}  is
that $\delta_{C_a}$ reduces to $\delta_C$ and $\Delta_{C_a}$
reduces to $\Delta_C$ under the natural projection $\pi: C_a \to
C_a/a_+C_a \cong C$. Assume that $C_a=C \otimes a$ as graded
coalgebras over $a$. We have a vector spaces decomposition
\[C_a= C \otimes a = C \otimes(k \oplus a_+) \\
=(C \otimes k)\oplus( C \otimes a_+) \\
=C \oplus(C \otimes a_+). \] Since $\delta_{C_a}$ reduces to
$\delta_C$ under projection $\pi$ we must have that
$\delta_{C_a}=\delta_C+e,$
where $e\in Coder(C \otimes a_+)$ is a degree one coderivation.
The fact that $\delta_{C_a}^N=0$ implies that $e$ should
satisfy certain identities which we call the $(N,M)$ Maurer-Cartan equation.

\bigskip

Let us review  the construction of $(N,M)$ Maurer-Cartan equation
\cite{AnDi1}. For $s=(s_1,...,s_n)\in \N^n$ we set
$l(s)=n$ and $|s|=\sum_i{s_i}$. For $1\leq i <n,$ let
$s_{>i}=(s_{i+1},...,s_n)$, for $1<i\leq n,$ let
$s_{<i}=(s_1,...,s_{i-1})$, also set $s_{>n}=s_{<1}=\emptyset$.
$\N^{(\infty)}$ denotes the set $\bigsqcup_{n=0}^{\infty}\N^n$,
where by convention $\N^{0}=\{\emptyset\}$. For $e\in Coder(C)$ we
set $e^{(s)}=e^{(s_1)}...e^{(s_n)}$, where $e^{(a)}=d^{a}(e)$ if
$a\geq 1$, $e^{(0)}=e$ and $e^{\emptyset}=1$. For $M\in\N$ we let
$$E_M=\{s\in\N^{(\infty)}:|s|+l(s)\leq M\}$$ and for $s\in E_M$ we
define integer $M(s)$ by $M(s)=M-|s|-l(s)$.\\

Recall that a discrete quantum mechanical
system is given by a directed graph together with a weight
attached to each of its edges. Let us introduce a discrete
quantum mechanical system  given by the following data
\begin{enumerate}
\item  The set of vertices is $\N^{(\infty)}$.
\item There is a unique directed edge
from vertex $s$ to vertex $t$ if and only if
$$t\in\{(0,s)\} \cup \{ s \} \cup \{ s+e_i\mid 1 \leq i \leq l(s) \}$$ where
$e_i=(0,..,\underset{\scriptsize{i-th}}{\underbrace{1}},..,0)\in\N^{l(s)}$.

\item An edge $e$ with source $s(e)$ and target $t(e)$  is weighted according to the table
\smallskip
\begin{center}
\begin{tabular}{|l|l|l|}\hline $s(e)$    &      $t(e)$ &
$v(e_i)$      \\  \hline $s$         &      $(0,s)$      & $1$
\\ \hline
$s$         &      $s$          & $(-1)^{|s|+l(s)}$    \\ \hline
$s$         &     $s+e_i$    & $(-1)^{|s_{<i}|+i-1}$ \\
\hline
\end{tabular}
\end{center}
\end{enumerate}

The set $P_M(\emptyset,s)$ consists of all paths
$\gamma=(e_1,...,e_M)$ such that $s(e_1)=\emptyset$ and
$t(e_M)=s$. The weight $v(\gamma)$ of $\gamma\in P_M(\emptyset,s)$
is $v(\gamma)=\prod_{i=1}^{M}v(e_i).$ \\

The $(N,M)$ Maurer-Cartan equation is $\sum_{k=1}^{M-1}c_k\delta^k=0,$ where
\[c_k=\sum_{\begin{subarray}{c} s\in E_M\\ M(s)=k \\ s_i<N\\ \end{subarray}}c(s,M)e^{(s)}
\hspace{.5cm}\text{and}\hspace{.1cm} c(s,M)=\sum_{\gamma\in
P_M(\emptyset,s)}v(\gamma).\]

In Section \ref{infinito}  we shall encounter the
problem of understanding the deformations of a codifferential.

\section{Differential graded algebras of depth N}

In this section we introduce the notion of differential graded algebras
of depth $N$ and construct explicitly the operad controlling such algebraic
structures.\\

Let $A$ be a $\Z$-graded vector space. Let $T(A[1])$ be the cofree
coalgebra cogenerated by $A[1]$, that is,
$$T(A[1])=\bigoplus_{n\geq 1}A[1]^{\otimes n}$$ with
$\Delta:T(A[1])\to T(A[1])$ given for $n \geq 2$ by
\[\Delta(a_1,...,a_n)=\sum (a_1,...,a_k)\otimes(a_{k+1},...,a_n).\]
We also consider  sub-coalgebras
\[T^{\leq k}(A[1])=\bigoplus_{1\leq n\leq k}A[1]^{\otimes n}.\]
Recall that a sequence of homomorphisms $m_k: A[1]^{\otimes k}\to
A[1]$ of degree one, for $k \in  \mathbb{N}_+$, defines a unique coderivation
$\delta:T(A[1])\to T(A[1])$ on $A[1]^{\otimes n}$ given by
\[\delta(a_1,....,a_n)=\sum_{k=1}^{n} \sum_{i=1}^{n-k+1}(-1)^{\sum_{j=1}^{i-1} \overline{a_i}} (a_1,...,a_{i-1},m_k(a_i,...,a_{i+k-1}),a_{i+k+1},...,a_n).\]

One may think of a dga-structure on a graded vector space $A$ as being a
codifferential $\delta\in Coder(T^{\leq 2}A[1])$ of degree one
given by $\delta=d + m.$ The condition $\delta^2=0$ is equivlent to the properties defining a differential
graded algebra.
Explicitly,
\begin{eqnarray*}
\delta^2&=&(d +m)(d + m)\\
&=&d^2+ dm + md+ m^2.
\end{eqnarray*}
We see that $\delta^2=0$ if and only if  $d^2=0$, $dm+md=0$ and $m^2=0$, i.e.,
$d$ is a square free map, satisfying the graded Leibnitz rule with respect to the associative product $m$.\\

We come to the main idea of this paper. Instead of looking at
codifferentials on $T^{\leq 2}(A[1])$, we may as well consider
$3$-codifferentials on $T^{\leq 2}(A[1])$, or more generally a
$N$-codifferentials  for $N \geq 2$. This idea leads to a new
mathematical entity which we define below.

\begin{defn}
Let $A$ be a $\Z$-graded vector space. A structure of graded
algebra of depth $N$ on $A$ is given by a degree one coderivation
$\delta\in Coder(T^{\leq 2}(A[1]))$ such that $\delta^N=0$.
\end{defn}

By the previous remarks a graded algebra of depth two is the same
as a differential graded algebra. From the decomposition of
$\delta=d + m$ in homogeneous components with $d:A[1] \to A[1]$ and $m:A[1]^{\otimes 2}\to
A[1]$, we can see that for a differential graded algebra of depth
$3$ we have
\begin{eqnarray*}
\delta^3&=&(d + m)(d + m)(d + m)\\
&=&(d+m)(d^2+dm+md+mm)\\
&=&d^3+d^2m+dmd+dmm+md^2\\
&+& mdm+mmd+mmm\\
&=&0.
\end{eqnarray*}
Looking at the homogeneous components of the identity above, we conclude that $d$ and
$m$ should satisfy:
\begin{enumerate}
\item  $m^3=0$. It is natural to call this property the condition
of $3$-associativity for $m$.
\item $d$ satisfies two generalized graded Leibnitz rules
\[d^2m+dmd+md^2=0,\]
\[dmm+mdm+mmd=0.\]
\item $d^3=0$, i.e., $d$ is a $3$-differential on $A$.
\end{enumerate}

Before we continue with our study of algebras of depth $N$ we
introduce a few combinatorial tools.

\begin{defn}
A finite planar rooted tree $\Gamma$ consists of a pair of finite
sets $V_\Gamma$ and $E_\Gamma$, called the vertices and edges of
$\Gamma$, and maps $s,t:E_\Gamma\to V_\Gamma$ satisfying:
\begin{enumerate}
\item  There exists a unique distinguished vertex
$r \in V_\Gamma$ called the root of $\Gamma$, such that
$|t^{-1}(v_r)|=1$ and $|s^{-1}(v_r)|=0$.\footnote{ $|X|$ denotes
the cardinality of a finite set $X$.}
\item There is a unique path from each vertex to the root.
\end{enumerate}
A vertex $v \in V_\Gamma$ such that $|t^{-1}(v)|=0$ is called a
leave.  A vertex $v$ the is neither the root nor a leave is called internal.
A vertex $v$ is said to be $n$-ary  if
$|t^{-1}(v)|=n$. $l(T)$ is the set of leaves in $T$.
\end{defn}

We use planar rooted trees to encode in a simpler notation operators on $A$.
First we associate rooted trees with the operators $d$ and $m$ as done in
Section \ref{ndgla}. It should be clear that from these simple rules we can associate
with each tree $T$ whose internal vertices are either unary or
binary an operator $O_T$ from $A^{\otimes|l(T)|}$ to $A.$ Formally,
$O_T$ is defined inductively as follows
\begin{enumerate}
\item The operator associated with the unique rooted tree with two
vertices is the identity.

\item For any other tree $T$ let $v$ be the unique vertex connected
with one edge to the root. Assume that $v$ is unary. Then a tree $S$
is attached to $v$ and $O_T=d \circ O_S$.

\item Assume now that $v$ is binary. Then trees  $T_1$ and
$T_2$ are attached to $v$ and $O_T=m \circ (O_{T_1} \otimes
O_{T_2}).$
\end{enumerate}

The condition for $A$ to be a differential graded
algebra of depth $3$ may be described graphically as follows. The
condition $d^3=0$ corresponds with
\begin{center}
\includegraphics[height=2cm]{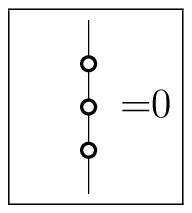}
\end{center}
The generalized Leibniz rules for $d$ and $m$ are
\begin{center}
\includegraphics[height=2cm]{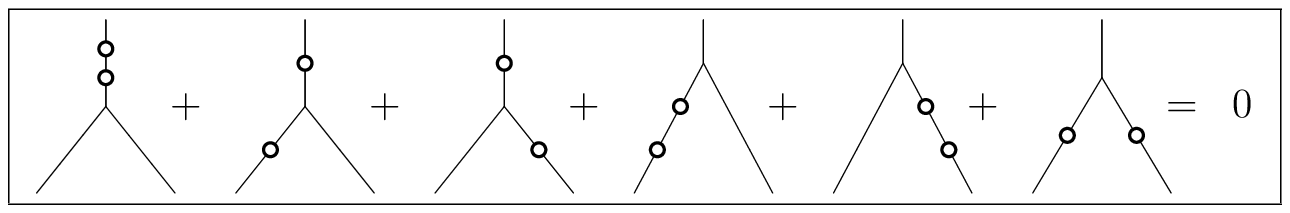}
\end{center}
and
\begin{center}
\includegraphics[height=4cm]{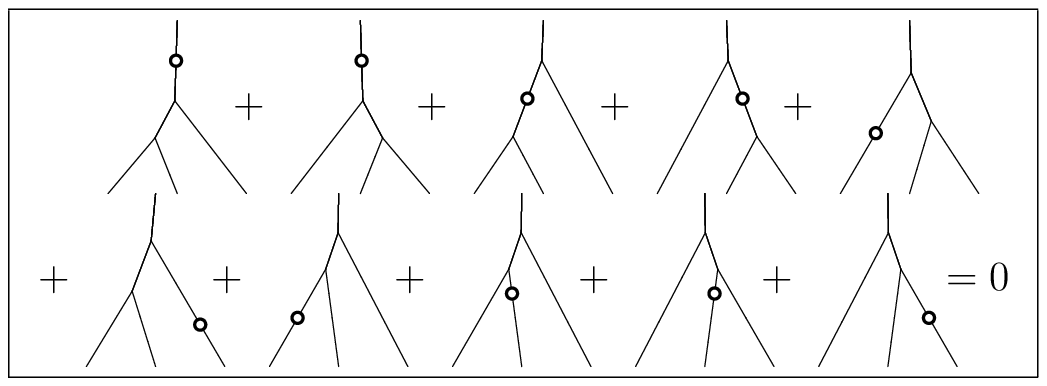}
\end{center}

The condition for $m$ to be $3$-associative is depicted as

\begin{center}
\includegraphics[height=2cm]{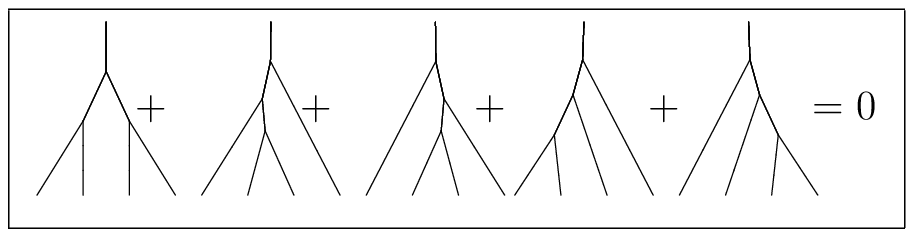}
\end{center}

\smallskip
Let us consider the graphical description of the axioms satisfied
by the operators defining a differential graded algebra of depth
$N$.
\begin{defn}
For $u+b=N$, Let $RT_l^{u,b}$ be the set of isomorphisms classes
of planar rooted tress with $l$ leaves, $u$ unary internal
vertices and $b$ binary internal vertices.
\end{defn}

We are ready to give a graphical description of the equations
defining a differential graded algebra of depth $N$.

\begin{thm}\label{graph}
Maps $d:A[1] \to A[1]$ and $m_A:A[1]^{\otimes 2}\to A[1]$ of
degree one define a differential graded algebra of depth $N$
structure on $A$ if and only if for $l=1,...,N+1$ and $u+b=N$ the
following identities hold
$$\sum_{T
\in RT_l^{u,b}}O_T=0.$$
\end{thm}
\begin{proof} Follows from the fact that
if $T \in RT_l^{u,b}$ then $O_T$ is a degree $u$ operator, and from
the identity $$(d + m)^N=\sum_{l=1}^{N+1}\sum_{u+b=N} (
\sum_{T
\in RT_l^{u,b}}O_T).$$
\end{proof}

Theorem \ref{graph} has several interesting consequences. First,
it allows us to define an operad over $gvect$  whose algebras in
$gvect$ are differential graded algebras of depth $N$. Second, the
existence of such an operad implies the existence of the free
differential graded algebra of depth $N$ generated by a graded
vector space.

\begin{defn}
The operad $\mbox{dgass}^N$ is given  by $
\mbox{dgass}^N= \overline{\mbox{dgass}}^N/I,$ where
$$\overline{\mbox{dgass}}^N(n)= <(T,f) \mid T \in RUBT_n \mbox{ and } f:l(T) \rightarrow [n] \mbox{ bijection}>,$$
and $RUBT_n=\sqcup_{u=0,b=0}^{\infty} RT_n^{u,b}.$ An element $(T,f)$ with $T \in
RT_n^{u,b}$ is placed in degree $u$. Compositions are given by
grafting of trees. $I$ is the operadic ideal $$I= <\sum_{T \in
RT_n^{u,b}}(T,f)
\mid u+b=N>.$$
\end{defn}

\begin{prop}
\begin{enumerate}
\item A $\mbox{dgass}^N$-algebra is a differential graded algebra of depth $N$.

\item $\bigoplus_{n=0}^{\infty}\mbox{dgass}^N(n) \otimes_{S_n}
V^{\otimes n}$ is the free differential graded algebra of depth
$N$ generated by the graded vector $V$.
\end{enumerate}
\end{prop}

\section{Deformations of N-associative algebras}

\medskip

In this section we study $N$-associative algebras. We define the operad
whose algebras are $N$-associative algebras and study the
infinitesimal deformations of such structures.

\begin{defn}
A $N$-associative algebra is a vector space $A$ together with a
linear map $m:A \otimes A \longrightarrow A$ such that $m^3=0.$
\end{defn}

Explicitly a product on a vector space $A$ is $3$-associative if
for $a,b,c,d
\in A$ we have that \[(ab)(cd)+ (a(bc))d + ((ab)c)d =  a((bc)d) + a(b(cd)).\]

Let us provide an elementary example of a $3$-associative
algebra.

\begin{exmp}
Let $A$ be the free non-associative algebra generated by
$a,b,c,d$ subject to the relations: $a^2=b, ab=d, ba=c$, the
product of two other letters is zero. $A$ is $3$-associative
algebra.
\end{exmp}

Let us now define explictly the operad on $vect$ whose algebras in $vect$ are
$N$ associative algebras. Let $RBT_n$ be the set of isomorphisms
classes rooted trees with binary internal vertices and $n$ leaves.

\begin{defn}
The operad $\mbox{ass}^N$ is given by $
\mbox{ass}^N= \overline{\mbox{ass}}^N/I,$ where
$$\overline{\mbox{ass}}^N(n)= <(T,f) \mid T \in RBT_n \mbox{ and } f:l(T) \rightarrow [n] \mbox{ bijection}>.$$
Compositions are given by grafting of trees. $I$ is the operadic
ideal $$I= <\sum_{T \in RBT_{N+1}}(T,f)>.$$
\end{defn}

\begin{prop}
\begin{enumerate}
\item An  $\mbox{ass}^N$-algebra is the same as a $N$-associative algebra.

\item $\bigoplus_{n=0}^{\infty}\mbox{ass}^N(n) \otimes_{S_n}
V^{\otimes n}$ is the free $N$-associative algebra generated by
the vector $V$.
\end{enumerate}
\end{prop}

\bigskip

It would be interesting to compute the generating series
of the operad $\mbox{ass}^N$. For  $2 \leq n \leq N$ the dimension of $\mbox{ass}^N(n)$ is $n!C_{n-1}$, where
$C_{n-1}=\frac{1}{n}\binom{2n-2}{n-1}$  is the Catalan number.  The dimension of $\mbox{ass}^N(N+1)$ is
$\frac{1}{N!}\binom{2N}{N}-(N+1)!.$\\

Notice that $ass^N$ is an $(N+1)$ homogeneous operad in the sense
that it is the quotient of a free operad $\overline{\mbox{ass}}^N$
generated by elements in $\overline{\mbox{ass}}^N(2)$ by an ideal
generated by elements in $\overline{\mbox{ass}}^N(N+1).$ Koszul duality
was originally introduced by Priddy \cite{P}
in the context of quadratic algebras. Ginzburg and Kapranov in
their seminal paper \cite{G} defined Koszul dualtity for quadratic operads.  Berger  in \cite{B1}
and Berger, Dubois-Violette
and Wambst in \cite{B2},  developed a theory of Koszul duality for
$N$-homogeneous algebras. It is natural to wonder if there exists
a notion of Koszul duality for $N$-homogeneous operads, and in
particular  what the Koszul dual of $ass^N$ might be.\\

Our next goal is to define the analogue for $N$-associative
algebras of the fact that infinitesimal deformations of an
associative algebra are controlled by its second Hochschild cohomology
group.

\begin{defn}
\begin{enumerate}
\item For a graded vector space $V$, let  $C(V,V)$ be the graded vector
space given by $C(V,V)=\bigoplus_{n=1}^{\infty}Hom(V^{\otimes
n},V)[1-n]$.
\item For $(A,m)$  a differential graded algebra of depth $N$  and  $k \geq N$, let
$$t_k:C^2(A,A) \longrightarrow C(A,A) \mbox{ be given by } t_k(f)= \sum_{i=0}^{k-1}m^i f m^{k-i}.$$
\item The cohomology group $H_{N,M}^2(A,A)$ is given by $$H_{N,M}^2(A,A)=
Ker(t_{M}:C^2(A,A) \longrightarrow C(A,A))/Im(t_1:C^1(A,A)
\longrightarrow C^2(A,A)),$$
where $t_1:C^1(A,A)
\longrightarrow C^2(A,A))$ is given by $t_1(g)=mg-gm.$
\end{enumerate}
\end{defn}

\smallskip
The definition above is consistent by the following result.

\begin{lem}\label{imp}
$t_{k} \circ t_1=0$ for $k \geq N.$
\end{lem}
\begin{proof}  For $g
\in C(A,A)$ we have
\begin{eqnarray*}
t_{k} \circ t_1(g)&=& t_{k}(t_1(g)) \\
&=&\sum_{i=0}^{k-1}m^i t_1(g) m^{k-1-i}\\
&=&\sum_{i=0}^{k-1}m^i (mg - gm )m^{k-1-i}\\
&=&\sum_{i=1}^{k}m^{i}g m^{k-i}- \sum_{i=0}^{k-1}m^i g m^{k-i}\\
&=&0.
\end{eqnarray*}
\end{proof}
\medskip

We are ready to state the main result of this section.\\

\begin{thm}\label{imp2}
\begin{enumerate}
\item Infinitesimal deformations of a $N$-associative algebra
$(A,m)$ into a $M$-associative algebra are determined by the cohomology group $H_{N,M}^2(A,A)$.
\item There are  inclusion maps $H_{N,M}^2(A,A) \longrightarrow H_{N,M+1}^2(A,A).$
\item Infinitesimal deformations of a $N$-associative algebra into a proper
$M$-associative algebra, $M > N$, are determined by the quotient space $H_{N,M}^2(A,A)/H_{N,M-1}^2(A,A).$

\end{enumerate}
\end{thm}
\begin{proof} Let $m_h(a,b)=m(a,b)+hf(a,b)$ be and infinetisimal deformation
of $A$ into a $M$-associative algebra with $f:A\otimes A\to A$ and $h$ a formal variable such
that $h^2=0$. The condition for $m_h$ to be $M$-associative is
that $m_h^M=0$. We have
$$m_h^M = (m+hf)^M
=h\sum_{i=0}^{M-1} m^i f m^{M-1-i}=h t_{M}(f).$$ Thus $m_h$ is
an infinitesimal deformation of $m$ into a $M$-associative product if and only if $f$ is a $t_M$ closed
element in $C^2(A,A),$ i.e., $f$ belongs to
$Ker(t_{M})$.\\

Any $g \in C^1(A,A)$ is a linear map $g:A\to A$ and gives rise to a
formal isomorphisms connected to the identity  $\rho=I+hg$
with inverse  $\rho^{-1}=I-hg$. Thus any $g
\in C^1(A,A)$  defines a infinitesimal deformation of $m$
given by
$m_g=\rho \circ m \circ ( \rho^{-1} \otimes \rho^{-1} ).$
We have that
\begin{eqnarray*}
m_g(a,b)&=&\rho(m(\rho^{-1}(a),\rho^{-1}(b)))\\
&=&\rho(m(a-hg(a),b-hg(b)))\\
&=&\rho(m(a,b))-h[\rho(m(g(a),b))+\rho(m(a,g(b)))]\\
&=&m(a,b)+h[g(m(a,b))-m(g(a),b)-(m(a,g(b)))]\\
&=&m(a,b)+ h[t_1(g)(a,b)]
\end{eqnarray*}
We see that $m_g$ differs from $m$ by an element $Im(t_1)$. We
have shown that $H_{N,M}^2(A,A)$ controls the infinitesimal
deformations of $m$ into a $M$-associative product. Properties $2$ and
$3$ follow from property $1$.
\end{proof}

Let us apply this result to the case of associative algebras.

\begin{cor}
An associative algebra admits an infinitesimal
deformation into a proper  $3$-associative algebra  if and only if
there exists a map $f:A
\otimes A
\longrightarrow A$ such that
\begin{enumerate}
\item For all $a,b,c,d \in A$ the following identity holds
$$f(ab,c)d +af(bc,d)+f(a,b)cd = abf(c,d)+f(a,bc)d +af(b,cd).$$

\item The following identity does not hold for all $a,b,c \in A$
$$f(a,b)c +f(ab,c) =af(b,c) + f(a,bc). $$
\end{enumerate}
\end{cor}
\begin{proof} $f$ defines a infinitesimal deformation
of $m$ into a $3$-associative product if $mfm=0$. It is also a infinitesimal deformation into an
associative product if $mf+fm=0$.
\end{proof}

\section{Introduction to
$A_{\infty}^{N}$-algebras}\label{infinito}

In this section we  introduce the notion of $A_{\infty}^N$-algebras which
generalizes both $N$-associative algebras and $A_{\infty}$
algebras.\ \ The interested reader may consult \cite{Fuk}, \cite{KaSt} , \cite{Ke} and \cite{So} for more
on $A_{\infty}$-algebras.
\medskip

\begin{defn}
A structure of $A_\infty^{N}$-algebra on a graded vector space $A$ is
given by a sequence of degree one maps $m_k: A[1]^{\otimes k}\to A[1]$
, for $k \in \mathbb{N}_+$, such that the associated
coderivation $\delta=\sum m_i$ on $T(A[1])$ satisfies
$\delta^{N}=0$.
\end{defn}

\medskip

The condition $\delta^{2}=0$ defining a $A_{\infty}^2$-algebra is usual condition for an $A_\infty$-algebra
\[\sum_{r+s+t=n}m_{r+1+t}\circ (1^{\otimes r}\otimes m_{s}\otimes 1^{\otimes t})=0.\]
The condition  $\delta^{3}=0$ defining an $A_{\infty}^3$-algebra
is

\begin{eqnarray*}
\sum_{a+b+c+d+e=n}m_{a+e+1}\circ(1^{\otimes a}\otimes m_{b+d+1}\otimes 1^{\otimes e})\circ
(1^{\otimes a+b}\otimes m_{c}
\otimes 1^{\otimes d+e})&+&\\
m_{a+b+c+e+1}\circ (1^{\otimes a}\otimes m_{b}\otimes 1^{\otimes c+e+1})\circ ( 1^{\otimes a+b+c} \otimes
m_{d}\otimes 1^{\otimes e})&+&\\
m_{a+c+d+e+1}\circ (1^{\otimes a+c+1}\otimes m_{d}\otimes 1^{\otimes e})\circ ( 1^{\otimes a} \otimes
m_{b}\otimes 1^{\otimes c+d+e})&=&0.
\end{eqnarray*}

It becomes difficult to write explicitly the condition
$\delta^{N}=0$ for $N\geq 4$, so we shall write it in terms of trees.
Denote by $RT_{l}^n$ the set of isomorphisms classes of rooted
planar trees with $l$ leaves and $n$ internal vertices. For
example the following trees are in $RT_{16}^3$.

\begin{center}
\includegraphics{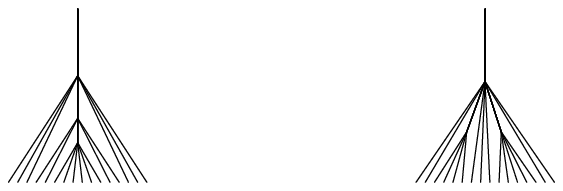}
\end{center}

The condition $\delta^{N}=0$ holds if and only if for each $l \in \mathbb{N}_+$
\[\sum_{\Gamma\in RT_{l}^N}m_{\Gamma}=0,\]
where $m_{\Gamma}$ is defined by a procedure similar to that
explained in the previous section, i.e.,  putting the $m_s$
operator on each vertex with $s$ incoming edges attached to it.\\

\begin{defn}
\begin{enumerate}
\item An $A^{N}_{\infty}$-morphism between  $A^{N}_{\infty}$-algebras
$(A,\delta_A)$ and $(B,\delta_B)$ is a degree zero coalgebra morphism
$f:T(A[1])\to T(B[1])$ such that $f\circ
\delta_A=\delta_B\circ f$.

\item Let $(A,\delta)$ be $A^{N}_{\infty}$-algebra. Let $m_1$ is the $Hom(A[1],A[1])$ component of $\delta$.
It is not hard to check that $m_1^N=0$. The
cohomology of $A$ is the cohomology of the $N$-complex $(A,m_1)$.

\item $A^{N}_{\infty}$-algebras $A$ and $B$
are quasi-isomorphic if there exist an $A^{N}_{\infty}$-morphism
between them inducing an isomorphism in cohomology.
\end{enumerate}
\end{defn}

The following result gives examples of $A_{\infty}^N$-algebras.

\begin{lem}
Let $(A,m,d)$ be a graded algebra of depth $N$. Setting
$m_1(a)=d(a)$, $m_2(a,b)=ab$ and $m_k=0$ for $k\geq 3$  defines a
$A^N_{\infty}$-structure on $A$.
\end{lem}

One can define the category of $A_{\infty}^{Nil}$-algebras
whose objects are $A_{\infty}^N$-algebras for some $N$. Next
result shows that the category of $A_{\infty}^{Nil}$-algebras is
monoidal.

\begin{thm}
Let $(A,\delta_A)$  be an $A^{N}_{\infty}$-algebra and $(B,\delta_B)$
be an $A^{M}_{\infty}$-algebras. Setting
$\delta_{\scriptsize{A\otimes B}}=\delta_A\otimes Id+Id\otimes
\delta_B$, we get that $(A\otimes B,\delta_{\scriptsize{A\otimes B}})$ is
a $A^{(N+M-1)}_{\infty}$-algebra.
\end{thm}

Let us define the graded operad $A_{\infty}^{N}$
whose algebras in $gvect$ are $A_{\infty}^{N}$-algebras.

\begin{defn}
The operad $A_{\infty}^N$ is given by $
A_{\infty}^N= \overline{A_{\infty}}^N/I,$ where for $n \in \mathbb{N}_+$
$$\overline{A_{\infty}}^N(n)= <(T,f) \mid T \in RT_n \mbox{ and } f:l(T) \rightarrow [n] \mbox{ bijection}>.$$
If $T$ has $v_k$ vertices with valence $k$, then $(T,f)$ is placed in degree $\sum_{k}v_k(2-k)$.
Compositions are given by grafting of trees. $I$ is the operadic
ideal $$I= <\sum_{T \in RT_{l}^N}(T,f)>.$$
\end{defn}

\smallskip

\begin{prop}
\begin{enumerate}
\item Algebras over $A_{\infty}^N$ in $gvect$ are $A_{\infty}^N$-algebras.

\item $\bigoplus_{n=1}^{\infty}\mbox{ass}^N(n) \otimes_{S_n}
V^{\otimes n}$ is the free $A_{\infty}^N$-algebra generated by the graded space $V$.
\end{enumerate}
\end{prop}

\smallskip

Let us study infinitesimal deformations of $A_{\infty}^N$-algebras.

\begin{defn}
\begin{enumerate}

\item For $(A,\delta)$  an $A_{\infty}^N$-algebra  and  $k \geq N$, let
$$t_k:C(A,A) \longrightarrow C(A,A) \mbox{ be given by } t_k(f)= \sum_{i=0}^{k-1}m^i f m^{k-i}.$$
\item The cohomology group $H_{N,M}(A,A)$ is given by $$H_{N,M}(A,A)=
Ker(t_{M}:C(A,A) \longrightarrow C(A,A))/Im(t_1:C^1(A,A)
\longrightarrow C(A,A)),$$
where $t_1:C^1(A,A)
\longrightarrow C(A,A)$ is given by $t_1(g)=mg-gm.$
\end{enumerate}
\end{defn}

\smallskip
Next couple of results are proved as Lemma \ref{imp} and Theorem
\ref{imp2}, respectively.

\begin{lem}\label{imp*}
$t_{k} \circ t_1=0$ for $k \geq N.$
\end{lem}

\begin{thm}\label{imp2*}
\begin{enumerate}
\item Infinitesimal deformations of an $A_{\infty}^N$-algebra
$A$ into an $A_{\infty}^M$-algebra are determined by the cohomology group $H_{N,M}(A,A)$.
\item There are inclusion maps $H_{N,M}(A,A) \longrightarrow H_{N,M+1}(A,A).$
\item Infinitesimal deformations of an $A_{\infty}^N$-algebra into a proper
$A_{\infty}^M$-algebra, $M > N$, are determined by the quotient space $H_{N,M}(A,A)/H_{N,M-1}(A,A).$

\end{enumerate}
\end{thm}

Let us apply this result to  $A_{\infty}$-algebras.

\begin{cor}
An $A_{\infty}$-algebra $A$ admits an infinitesimal
deformation into a proper $A_{\infty}^M$-algebra  if and only if
there exists $f \in C(A,A)$ such that the following identity holds
\begin{eqnarray*}
\sum_{a+b+c+d+e=n}m_{a+e+1}\circ(1^{\otimes a}\otimes f_{b+d+1}\otimes 1^{\otimes e})\circ
(1^{\otimes a+b}\otimes m_{c}
\otimes 1^{\otimes d+e})&+&\\
m_{a+b+c+e+1}\circ (1^{\otimes a}\otimes f_{b}\otimes 1^{\otimes c+e+1})\circ ( 1^{\otimes a+b+c} \otimes
m_{d}\otimes 1^{\otimes e})&+&\\
m_{a+c+d+e+1}\circ (1^{\otimes a+c+1}\otimes m_{d}\otimes 1^{\otimes e})\circ ( 1^{\otimes a} \otimes
f_{b}\otimes 1^{\otimes c+d+e})&=&0.
\end{eqnarray*}

but the following condition fails
\[\sum_{a+b+c=n}m_{a+1+b}\circ (1^{\otimes a}\otimes f_b \otimes 1^{\otimes c}) +
\sum_{a+b+c=n}f_{a+1+b}\circ (1^{\otimes a}\otimes m_b \otimes 1^{\otimes c})=0.\]
\end{cor}

We consider full deformations of $A_{\infty}^N$-algebras. Let $k$ be a field and consider a local $k$-algebra $a$ such that
$k \cong a/a_+$ where $a_+$ is the unique maximal ideal in $a$. Then $a \cong a\oplus a_+$ as
vector spaces.

\begin{defn}
Let $(A,\delta)$ be an $A_{\infty}^N$-algebra. A deformation
of $(A,\delta)$ over $a$ is a $A_{\infty}^N$-algebra
$(A_a,\delta_a)$ over $a$ such that
$A_a/a_+A_a$ is isomorphic to $A$ as
$A_{\infty}^N$-algebra.  Equivalently,
$\delta_a$ maps to $\delta$ under the natural projection
$\pi: A_a \to A_a/a_+ A_a\cong A$.
\end{defn}

Suppose that $A_a=A\otimes a$ as $a$-modules. Then $\delta_a=\delta+e$ where $e\in
Coder(T(A[1])\otimes a_+)$ is a coderivation of degree one. The next result follows from
the discussion on the deformations of $N$-codifferntials at the end of Section \ref{co}.

\begin{thm}
Let $e \in Coder(T(A[1])\otimes a_+)$.  $\delta + e$ defines an $A_{\infty}^M$-algebra
structure on $A_a=A\otimes a$ if and only if the $(N,M)$ Maurer-Cartan equation
$\sum_{k=0}^{M-1}c_k\delta^k=0$ holds, where
\[c_k=\sum_{\begin{subarray}{c} s\in E_M\\ N(s)=k \\ s_i<N\\ \end{subarray}}c(s,M)e^{(s)}
\hspace{.5cm}\text{and}\hspace{.4cm} c(s,M)=\sum_{\gamma\in
P_M(\emptyset,s)}v(\gamma).\]
\end{thm}

We close this paper with a discussion of some open problems. The operadic description allows us to construct free $A^{N}_{\infty}$-algebras,
however further examples are needed. One can define a notion of $L^{N}_{\infty}$-algebras along
the lines of our definition of $A^{N}_{\infty}$-algebras. It would be interesting
to describe  explicitly $L^{N}_{\infty}$-algebras and find examples of such structures.
There are of course many open questions redarding $A^{N}_{\infty}$-algebras and $L^{N}_{\infty}$-algebras.
What are the generating series of $A^{N}_{\infty}$ and $L^{N}_{\infty}$?
Is any $A^{N}_{\infty}$-algebra quasi-isomorphic to a
differential graded algebra of depth $N$?
$A_{\infty}$-algebras can be defined geometrically via Stasheff's
associahedra \cite{Sta}. Is there an analogous geometrical description for
$A^{N}_{\infty}$-algebras? Notice that the equation defining $A^{N}_{\infty}$-algebras
describes the $N-1$ differential of the operation $m_k$ with $k \geq 2$ as linear sum of operators constructed
from the composition of $N$ operators $m_s$, each differentiated at most $N-2$ times.

\bibliography{bib}

\bigskip

\noindent ragadiaz@gmail.com \\
\noindent mangel@euler.ciens.ucv.ve  \\

\end{document}